\documentclass[11pt,leqno]{article} 
\usepackage{graphics}

\usepackage{lipsum}

\let\OLDthebibliography\thebibliography
\renewcommand\thebibliography[1]{
  \OLDthebibliography{#1}
  \setlength{\parskip}{2pt}
  \setlength{\itemsep}{2pt plus 0.3ex}
}

\newtheorem{thm}{Theorem}[section]

\newtheorem{prop}{Proposition}

\newcommand{\beqa}{\begin{eqnarray}}
\newcommand{\eeqa}{\end{eqnarray}}

\newcommand{\pf}{\noindent {\bf Proof:} $\s$ }
\newcommand{\epf}{ \hfill$\diamondsuit$ \medskip}

\newcommand{\ds}{\displaystyle}
\newcommand{\beq}{\begin{equation}}
\newcommand{\eeq}{\end{equation}}
\newcommand{\lbl}{\label}
\newcommand{\s}{\; \;}

\newcommand{\la}{\lambda}

\newcommand{\ra}{\rightarrow}
\newcommand{\al}{\alpha}
\newcommand{\p}{\varphi}

\title{Necessary and sufficient condition for existence at resonance for eigenvalues of multiplicity two }

\author{
Philip Korman   \\ 
Department of Mathematical Sciences \\ 
University of Cincinnati \\ 
Cincinnati Ohio 45221-0025 \\
}

\date{}

\begin{document}

\maketitle
\begin{abstract} 
We establish necessary and sufficient condition for existence of solutions for a class of semilinear Dirichlet problems with the linear part  at resonance  at eigenvalues of multiplicity two. The result is applied to give a condition for  unboundness of all solutions of the corresponding semilinear heat equation.
 \end{abstract}

\begin{flushleft}
\end{flushleft}

\begin{flushleft}
\end{flushleft}

\section{Introduction}
\setcounter{equation}{0}
\setcounter{thm}{0}
\setcounter{lma}{0}

The study of elliptic problems at resonance was initiated by the classical paper of E.M. Landesman and A.C. Lazer \cite{L}. On a bounded smooth domain $D \subset R^n$ consider the Dirichlet problem 
\beqa
\lbl{0}
& \Delta u+\la _k u+g(u)=f(x) \s \mbox{for $x \in D$} \,, \\ \nonumber
& u=0 \s \mbox{on $\partial D$} \,.
\eeqa
Here $\la _k$ is an eigenvalue of the Laplacian $\Delta$ on $D$ with zero boundary condition, so that the problem is {\em at resonance}. The function $f(x) \in L^2 (D)$ is  given.
For the nonlinear term $g(u)$ it is assumed that the limits $g(\infty)$ and $g(-\infty)$ exist and  
\beq
\lbl{1}
g(-\infty) <g(u)<g(\infty),  \s \mbox{for all $u \in (-\infty,\infty)$} \,.
\eeq
Let us recall the classical theorem of E.M. Landesman and A.C. Lazer \cite{L} in the form of S.A. Williams \cite{wi}.

\begin{thm}(\cite{L},\cite{wi}) \label{thm:1}
Assume that $g(u)$ satisfies (\ref{1}), $f(x) \in L^2(D)$, while for any $w(x)$
belonging to the eigenspace of $\lambda _k$
\begin{equation} \label{2}
\int _D f(x) w(x) \, dx<g(\infty) \int _{w>0} w \, dx
+g(-\infty) \int _{w<0} w \, dx \,.
\end{equation}
Then the problem (\ref{0}) has a solution $u(x) \in W^{2,2}(D) \cap W^{1,2}_0(D)$. Condition (\ref{2}) is also necessary
for the existence of solutions.
\end{thm}

Originally E.M. Landesman and A.C. Lazer \cite{L} assumed additionally that the eigenvalue $\la _k$ is simple. Soon, S.A. Williams \cite{wi} produced the more general statement given above. However,  for a while no examples for multiple eigenvalues were known, until we observed in \cite{K} that another classical result of A.C. Lazer and D.E. Leach \cite{L1} on periodic solutions of semilinear harmonic oscillator  provides an example to  Theorem \ref{thm:1} in case of double eigenvalues (giving incidentally another proof of  Lazer-Leach theorem, in addition to a number of other known proofs, see e.g., \cite{H} and \cite{K}).
We showed in \cite{K} that while the necessary condition of Lazer-Leach result is easy to prove, the sufficiency part follows by verifying  the condition (\ref{2}), and applying Theorem \ref{thm:1}.
\medskip

In this paper we prove a similar result for a  disc in $R^2$, thus providing the first PDE example for Theorem \ref{thm:1} in case of a multiple dimensional eigenspace.  Even for simple domains the eigenspace of a multiple eigenvalue can be very complicated, and multiplicity of eigenvalues may vary in non-obvious ways. So that verifying the inequality (\ref{2}) for any element $w(x)$ of the eigenspace appears to be very challenging for other domains (the integrals $\int _{w>0} w(x) \, dx$ and $\int _{w<0} w(x) \, dx$ are unlikely to remain constant over an eigenspace).
\medskip

\noindent
{\large \bf Example }  Let $D=(0,\pi) \times (0,\pi)$ in $R^2$. The eigenvalues of 
\[
\Delta u +\la u=0 \,, \s \mbox{in $D$}\, \s u=0 \s \mbox{on $\partial D$}
\]
are $\la _{nm}=n^2+m^2$ with positive integers $n$ and $m$, corresponding to the eigenfunctions $\sin nx \sin my$, see e.g., \cite{pr}. These eigenfunctions are obtained by separation of variables, and there are no other eigenfunctions since these eigenfunctions form a complete set in $L^2(D)$.   The principal eigenvalue $\la _1=2$ is simple, with the corresponding eigenfunction $\sin x \sin y>0$. The eigenvalue $\la _2=5=1^2+2^2$ has multiplicity two, with
the eigenspace spanned by $\sin x \sin 2y, \sin 2x \sin y$. The eigenvalue $\la _3=8=2^2+2^2$ is simple, with
the eigenspace spanned by $\sin 2x \sin 2y$. The eigenvalue $50=1^2+7^2=5^2+5^2$ is triple, with
the eigenspace spanned by $\sin x \sin 7y, \sin 7x \sin y, \sin 5x \sin 5y$. The eigenvalue $65=1^2+8^2=4^2+7^2$ is quadruple,  the eigenvalue $325=1^2+18^2=6^2+17^2=10^2+15^2$ has multiplicity six, and so on. It is natural to ask if there is an eigenvalue of any multiplicity. In number theoretic terms the possible conjecture is: for any even integer $2m$ one can find an integer $N$ that can be represented as $N=p^2+q^2$, with integers $p \ne q$, in exactly $m$ different ways, while  for any odd integer $2m+1$ one can find an integer $M$ that can be represented as $M=p^2+q^2$, with integers $p \ne q$, in exactly $m$ different ways, and in addition, $M=r^2+r^2$ for some integer $r>0$.
\medskip

By contrast, for a disc $B _a: x^2+y^2 <a^2$ in $R^2$, we show that all eigenvalues of the Laplacian have multiplicity two, except for the principal one (which is simple), and that the integrals $\int _{w>0} w(x,y) \, dxdy$ and $\int _{w<0} w(x,y) \, dxdy$ remain constant over the entire eigenspaces, and can be explicitly calculated. We present a necessary and sufficient condition for the existence of solutions of the problem (\ref{0}) on $B _a$, for this case of resonance at a double eigenvalue. We prove the necessity part directly, while sufficiency is derived by verifying the conditions of Theorem \ref{thm:1}. Our result can be seen as a PDE analog of the Lazer-Leach theorem. As an application, we give a condition for unboundness of all solutions of the corresponding semilinear heat equation.
\medskip

Radial solutions on balls in $R^n$ were studied extensively, see e.g., P. Korman \cite{K2} or T. Ouyang and J. Shi \cite{OS}, stimulated by  the classical theorem of B. Gidas, W.-M. Ni and L. Nirenberg \cite{GNN} which asserts that any positive solution of semilinear Dirichlet problem on a ball is necessarily radially symmetric. Our result suggests that ball may be a special domain even when studying sign-changing non-symmetric solutions.
\medskip

Previously, P. Korman and D.S. Schmidt \cite{ks} studied resonance at the principal eigenvalue on $B$. They constructed $g(u)$ for which the problem has infinitely many solutions {\em for any} $f(x,y) \in L^2(B _a)$.

\section{Resonance for a two-dimensional disc}
\setcounter{equation}{0}
\setcounter{thm}{0}
\setcounter{lma}{0}

Remarkably, the eigenvalues of Laplacian on a disc $B _a : x^2+y^2 <a^2$ in two dimensions all have multiplicity two, except for the principal eigenvalue, which is simple.
Recall (see e.g.  \cite{pr}, p. $251$) that the  eigenvalues of the Laplacian on $B _a$ with zero boundary condition are $\ds \la _{n,m}=\frac{\al _{n,m}^2}{a^2}$ ($n=0,1,2,\ldots$; $m=1,2,\ldots$), with the corresponding eigenfunctions 
\beq
\lbl{5--}
J_n \left( \frac{\al _{n,m}}{a} r \right) \left(A  \cos n \theta+B \sin  n \theta  \right) \,, 
\eeq
where $\al _{n,m}$ is the $m$-th root of $J_n(x)$, the $n$-th Bessel function of the first kind, $r=\sqrt{x^2+y^2}$ ($A$ and $B$ are arbitrary constants). There are no other eigenfunctions, since the ones given above form a complete set in $L^2(B _a)$. 
The principal eigenpair is $\la _1=\frac{\al _{0,1}^2}{a^2}  \approx \frac{5.78}{a^2}$, $\p _1(r)=J_0(\frac{\al _{0,1}}{a} r)$.
One calculates $\la _2=\frac{\al _{1,1}^2}{a^2}  \approx \frac{14.62}{a^2}$, with $\al _{1,1} \approx 3.83$, and $\p _2=J_1 \left( \frac{\al _{1,1}}{a} r \right) \left(A  \cos  \theta+B \sin   \theta  \right)$, and so on, see the Example below. The principal  eigenvalue is simple, while all other eigenvalues have multiplicity two, because any two Bessel functions with indices different by an integer do not have any roots in common, see G.N. Watson \cite{w}, p. 484 for the following result.

\begin{prop}
\lbl{p1}
For any integers $n \geq  0$ and $m \geq  1$, the functions $J_n(x)$ and $J_{n + m}(x)$ have no common zeros other than the one at $x = 0$. 
\end{prop}

This result was apparently once a long standing conjecture (published in 1866), known  in the 19-th century as Bourget's hypothesis (after a 19th-century French mathematician),  until it was proved in 1929 by C.L. Siegel, see \cite{w}, and a very informative Wikipedia article on the Bessel functions. The name ``hypothesis" suggests that it was used to prove other results. It immediately implies the following result that we need.

\begin{prop}
For the disc $B _a$, all eigenvalues, other than the principal one, have multiplicity two.
\end{prop}

\pf
By Proposition \ref{p1}, all $\al _{n,m}$ are different, and hence the eigenspace of $\la _{n,m}$ is two-dimensional, and is given by (\ref{5--}).
\epf

It turns out that for any eigenvalue $\la _k$, $k \geq 2$, both integrals $\int _{w>0} w \, dx$ and $\int _{w<0} w \, dx$ on $B _a$ (appearing in (\ref{2})) remain constant for all $w(x)$ in the eigenspace of $\la _k$, and both integrals  can be easily calculated. Let $P_{n,m}$ denote the subset of $(0,a)$ where $J_n \left( \frac{\al _{n,m}}{a} r \right)>0$, and $N_{n,m}$ the subset of $(0,a)$ where $J_n \left( \frac{\al _{n,m}}{a} r \right)<0$. The following quantity
\beq
\lbl{5}
J_{n,m}=2 \int_{P_{n,m}} J_n \left( \frac{\al _{n,m}}{a} r \right) r \, dr-2 \int_{N_{n,m}} J_n \left( \frac{\al _{n,m}}{a} r \right) r \, dr \,.
\eeq
can be easily calculated using {\em Mathematica} for each pair of $n$ and $m$.

\begin{prop}\lbl{prop:3}
Let $w=J_n \left( \al _{n,m} r \right) \left(A  \cos n \theta+B \sin  n \theta  \right)$ be any element of the eigenspace of the eigenvalue $\ds \frac{\al _{n,m}^2}{a^2}>\la _1$, normalized so that $A^2+B^2=1$. Then on $B _a$
\[
\int_{w>0} w(r,\theta) \,  r  dr \, d \theta=J_{n,m} \,,
\]
\[
\int_{w<0} w(r,\theta) \,  r  dr \, d \theta=-J_{n,m} \,.
\]
\end{prop}

\pf
Write
\[
A  \cos n \theta+B \sin  n \theta=\sqrt{A^2+B^2}\cos \left(n \theta-\delta \right)= \cos \left(n \theta-\delta \right) \,,
\]
for some $\delta$. Then 
\beq
\lbl{6}
w=J_n \left( \al _{n,m} r \right) \cos \left(n \theta-\delta \right) \,.
\eeq
Let $P$ denote the set of $\theta \in (0,2\pi)$ where $\cos \left(n \theta-\delta \right)>0$, and $N$ the set of $\theta \in (0,2\pi)$ where $\cos \left(n \theta-\delta \right)<0$. It is easy to show that
\beqa
\lbl{6b}
& \int_P \cos \left(n \theta-\delta \right) \, d \theta=2 \,, \\ 
& \int_N \cos \left(n \theta-\delta \right) \, d \theta=-2 \,. \nonumber
\eeqa
Then, in view of (\ref{6}), (\ref{5})  and (\ref{6b})
\beqa \nonumber
& \int_{w>0} w(r,\theta) \,  r  dr \, d \theta= \int_{P_{n,m}\times P} w(r,\theta) \,  r  dr \, d \theta+\int_{N_{n,m} \times N} w(r,\theta) \,  r  dr \, d \theta \\ \nonumber
&=\int_{P_{n,m}} J_n \left( \frac{\al _{n,m}}{a} r \right) r \, dr\int_P \cos \left(n \theta-\delta \right) \, d \theta \\ \nonumber
& +\int_{N_{n,m}} J_n \left( \frac{\al _{n,m}}{a} r \right) r \, dr\int_N \cos \left(n \theta-\delta \right) \, d \theta=J_{n,m} \,, \nonumber
\eeqa
and similarly
\[
\int_{w<0} w(r,\theta) \,  r  dr \, d \theta= \int_{P_{n,m} \times N} w(r,\theta) \,  r  dr \, d \theta+\int_{N_{n,m} \times P} w(r,\theta) \,  r  dr \, d \theta=-J_{n,m} \,,
\]
completing the proof.
\epf

We now consider the problem (here $u=u(x,y)$)
\beqa
\lbl{8}
& \Delta u+\la _k u+g(u)=f(x,y) \,, \s \s \mbox{for $(x,y) \in B _a$} \,, \\ \nonumber
& u=0 \s \mbox{on $\partial B _a$} \,,
\eeqa
with the eigenvalue $\ds \la _k=\frac{\al ^2 _{n,m}}{a^2}>\la _1$ for some $n$ and $m$, corresponding to the eigenspace $J_n \left( \frac{\al _{n,m}}{a} r \right) \left(A  \cos n \theta+B \sin  n \theta  \right)$. Let us denote 
\beq
\lbl{7a}
\p _k=J_n \left( \frac{\al _{n,m}}{a} r \right)   \cos n \theta \,,
\eeq
\[
\psi _k=J_n \left( \frac{\al _{n,m}}{a} r \right)   \sin n \theta \,,
\]
\[
A_k=\int_{B _a} f(x,y) \p _k \, dxdy \,,
\]
\[
B_k=\int_{B _a} f(x,y) \psi _k \, dxdy \,.
\]
The numbers $A_k$ and $B_k$ can be easily approximated by {\em Mathematica} for any $f(x,y)$ and $k$.

\begin{thm}\lbl{thm:2}
Assume that $g(u)$ satisfies the condition (\ref{1}). Then the condition
\beq
\lbl{9}
\sqrt{A^2_k+B^2_k}<J_{n,m} \left(g(\infty)-g(-\infty) \right)
\eeq
is both necessary and sufficient for the existence of solution $u(x) \in W^{2,2}(B _a) \cap W^{1,2}_0(B _a)$ of (\ref{8}).
\end{thm}

\pf
(i) Necessity. Multiply (\ref{8}) by $\p _k$ and $\psi _k$ respectively and integrate
\beqa
\lbl{9a}
& A_k=\int _{B _a} g(u) \p _k\, dxdy \,,\\
& B_k=\int _{B _a} g(u) \psi _k\, dxdy \,. \nonumber
\eeqa
Multiply the first equation in (\ref{9a}) by $\frac{A_k}{\sqrt{A^2_k+B^2_k}}$, the second equation by $\frac{B_k}{\sqrt{A^2_k+B^2_k}}$,  and add the results. Denoting 
\[
w_k=J_n \left( \frac{\al _{n,m}}{a} r \right) \left(\frac{A_k}{\sqrt{A^2_k+B^2_k}}  \cos n \theta+\frac{B_k}{\sqrt{A^2_k+B^2_k}} \sin  n \theta  \right) \,,
\]
and using Proposition \ref{prop:3}, obtain
\[
\sqrt{A^2_k+B^2_k}=\int _B g(u) w _k\, dxdy<g(\infty) \int _{w_k>0}  w _k\, dxdy+g(-\infty) \int _{w_k<0}  w _k\, dxdy
\]
\[
=J_{n,m} \left(g(\infty)-g(-\infty) \right) \,.
\]

\noindent
(ii) Sufficiency. Assuming that (\ref{9}) holds, we shall verify the condition (\ref{2}) of Theorem \ref{thm:1}. Assuming that 
$\la _k=\frac{\al^2 _{n,m}}{a^2}$, let 
\[
w(x,y)=J_n \left( \frac{\al _{n,m}}{a} r \right) \left(A  \cos n \theta+B \sin  n \theta  \right)
\]
be any element of its eigenspace. By scaling $w$ in  (\ref{2}), we may assume that 
\[
A^2+B^2=1 \,.
\]
In view of Proposition \ref{prop:3} and (\ref{7a}), the condition (\ref{2}) of Theorem \ref{thm:1} that we need to verify  takes the form
\beqa \nonumber
& \int_{B _a} f(x,y) w(x,y) \, dxdy= AA_k+BB_k \\ \nonumber
& <J_{n,m} \left(g(\infty)-g(-\infty) \right) =g(\infty) \int _{w>0} w \, dxdy
+g(-\infty) \int _{w<0} w \, dxdy \,. \nonumber
\eeqa
Since $AA_k+BB_k \leq \sqrt{A^2_k+B^2_k}$, this inequality holds by (\ref{9}). By Theorem \ref{thm:1} the problem (\ref{8}) has a solution.
\epf

\noindent
{\bf Example } Consider the unit disc $x^2+y^2<1$, $a=1$. {\em Mathematica} readily returns zeroes of the Bessel functions
\beqa \nonumber
& \al _{0,1} \approx 2.40483 \,,\s  \al _{0,2} \approx 5.52008\,, \s \al _{0,3} \approx 8.65373 \,, \ldots \\ \nonumber
& \al _{1,1} \approx 3.83171\,, \s \al _{1,2} \approx 7.01559\,, \ldots  \\ \nonumber
& \al _{2,1} \approx 5.13562\,, \s \al _{2,2} \approx 8.41724\,, \ldots  \\ \nonumber
& \al _{3,1} \approx 6.38016\,, \s \al _{3,2} \approx 9.76102\,, \ldots \,. \\ \nonumber
\eeqa
The eigenvalues are $\la _1=\al^2 _{0,1}$, $\la _2=\al^2 _{1,1}$, $\la _3=\al^2 _{2,1}$, $\la _4=\al^2 _{0,2}$, $\la _5=\al^2 _{3,1}$, $\la _6=\al^2 _{1,2}$, and so on. Let us consider a case of resonance at  the sixth eigenvalue
\beqa
\lbl{18}
& \Delta u+\la _6 u+\frac{u}{\sqrt{u^2+1}}=f(x,y) \,, \s \s \mbox{for $x^2+y^2<1$} \,, \\ \nonumber
& u=0 , \s \mbox{on $x^2+y^2=1$} \,.
\eeqa
By above, the eigenspace of $\la _6$ is $ J_1 \left( \al _{_{1,2}} r \right) \left(A  \cos  \theta+B \sin   \theta  \right)$, with arbitrary numbers $A$ and $B$. The graph of $  J_1 \left(  \al _{_{1,2}} r \right)$ on $(0,1)$  has one root $r_0=\frac{\al _{1,1}}{\al _{1,2}}$, and it is positive on $P_{1,2}=(0,r_0)$, and negative on  $N_{1,2}=(r_0,1)$, see Figure \ref{fig:1}. By (\ref{5}), using {\em Mathematica}
\[
J_{1,2}=2 \int_0^{r_0} J_1 \left( \al _{_{1,2}} r \right) r \, dr-2 \int_{r_0}^1 J_n \left( \al _{_{1,2}} r \right) r  \, dr \approx 0.260759\,.
\]
For any $f(x,y)$, {\em Mathematica} can also easily compute  highly accurate approximation of the integrals
\[
A_6=\int _{x^2+y^2 <1} f(x,y)J_1 \left( \frac{\al _{1,2}}{a} r \right) \cos \theta \,dxdy \,,  
\]
\[
B_6=\int _{x^2+y^2 <1} f(x,y)J_1 \left( \frac{\al _{1,2}}{a} r \right) \sin \theta \,dxdy\,. 
\]
(Here $x=r \cos \theta$, $y=r \sin \theta$, and $dxdy=r \, dr d \theta$.)
Since $g(\infty)=1$ and $g(-\infty)=-1$, it follows by Theorem \ref{thm:2} that the problem (\ref{18}) has a solution if and only if 
\[
\sqrt{A^2_6+B^2_6}<2 J_{1,2} \,.
\]

\begin{figure}
\begin{center}
\scalebox{0.6}{\includegraphics{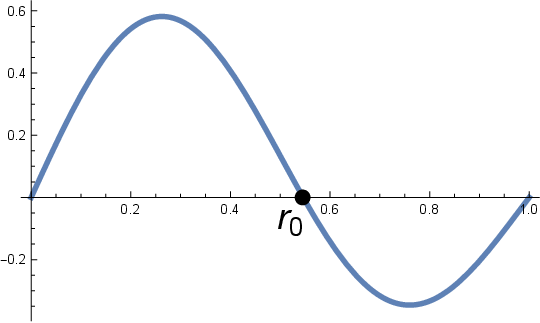}}
\caption{  The graph of $  J_1 \left( \al _{_{1,2}} r \right)$ on the interval $(0,1)$}
\lbl{fig:1}
\end{center}
\end{figure}

We now consider an application to the semilinear heat equation on a disc $B _a : x^2+y^2 <a^2$ (here $u=u(x,y,t)$)
\beqa
\lbl{11}
& u_t=\Delta u+\la _k u+g(u)-f(x,y) \,, \s \s \mbox{for $(x,y) \in B _a$}\,, t>0 \\\nonumber
&  u(x,y,t)=0 \,, \s\s \mbox{for $(x,y)$ on $\partial B _a$} \,, t>0\\ \nonumber
& u(x,y,0)=u_0(x,y)  \,, \ \nonumber
\eeqa
with given functions $f(x,y)$ and $u_0(x,y)$, and $g(u)$ satisfying (\ref{1}). Here $\la _k$, $k \geq 2$, is a double eigenvalue of Laplacian, as above. Steady states for this equation satisfy the equation (\ref{8}). By Theorem \ref{thm:2}, no steady states exist if 
\beq
\lbl{12}
\sqrt{A^2_k+B^2_k}>J_{n,m} \left(g(\infty)-g(-\infty) \right) \,.
\eeq
 
Denote $w_k=\frac{A_k}{\sqrt{A^2_k+B^2_k}} \p _k+\frac{B_k}{\sqrt{A^2_k+B^2_k}} \psi _k$, as above. (Recall that $
A_k=\int_{B _a} f(x,y) \p _k \, dxdy $,
$
B_k=\int_{B _a} f(x,y) \psi _k \, dxdy $.)

\begin{thm}\lbl{thm:3}
Assume that $g(u)$ satisfies the condition (\ref{1}), and that (\ref{12}) holds. Then solution of (\ref{11}) is unbounded for any initial data $u_0(x,y)$. In fact, defining $H(t)=\int_{B _a} u(x,y,t) w_k \, dxdy$, one has $H(t) \ra -\infty$ as $ t \ra \infty$.
\end{thm}

\pf
Multiply (\ref{11})  by $\frac{A_k}{\sqrt{A^2_k+B^2_k}} \p _k$ and integrate both sides over $B _a$
\[
\frac{A_k}{\sqrt{A^2_k+B^2_k}} \int _{B _a} u_t \p_k \, dxdy =\frac{A_k}{\sqrt{A^2_k+B^2_k}} \int _{B _a} g(u) \p_k \, dxdy-\frac{A^2_k}{\sqrt{A^2_k+B^2_k}}\,.
\]
Multiply (\ref{11})  by $\frac{A_k}{\sqrt{A^2_k+B^2_k}} \psi _k$, and integrate  over $B_a$
\[
\frac{A_k}{\sqrt{A^2_k+B^2_k}} \int _{B _a} u_t \psi_k \, dxdy =\frac{A_k}{\sqrt{A^2_k+B^2_k}} \int _{B _a} g(u) \psi_k \, dxdy-\frac{B^2_k}{\sqrt{A^2_k+B^2_k}}\,.
\]
Add the results, to obtain
\[
H'(t)=\int _{B _a} g(u) w_k \, dxdy- \sqrt{A^2_k+B^2_k}<-\epsilon \,,
\]
for some $\epsilon >0$, by  estimating the integral $\int _{B _a} g(u) w_k \, dxdy$  as in part (i) of Theorem \ref{1}, and using (\ref{12}). Then $H(t)<H(0)-\epsilon t$, concluding the proof.
\epf

In the ODE context related results on unbounded solutions were given by G. Seifert \cite{s}, J.M. Alonso and R. Ortega  \cite{O}, and P. Korman \cite{k1}.


\end{document}